\documentclass{article}

\usepackage[utf8]{inputenc} % allow utf-8 input
\usepackage[T1]{fontenc}    % use 8-bit T1 fonts
\usepackage{hyperref}       % hyperlinks
\usepackage{url}            % simple URL typesetting
\usepackage{booktabs}       % professional-quality tables
\usepackage{subcaption}
\usepackage{graphicx}
\usepackage{amsfonts}       % blackboard math symbols
\usepackage{nicefrac}       % compact symbols for 1/2, etc.
\usepackage{microtype}      % microtypography
\usepackage{xcolor}         % colors
\usepackage{multicol}

\usepackage{pgfplots}
\pgfplotsset{compat=1.16} 
\usepackage{amsmath}
\usepackage{amsthm}
\usepackage{amssymb}
\usepackage{mathrsfs}
\usepackage{textcomp}
\usepackage{authblk}
\usepackage{todonotes}
\usepackage{graphicx}
\usepackage{enumitem}

\numberwithin{equation}{section}

\theoremstyle{definition}
\newtheorem{defn}{Definition}[section]

\theoremstyle{remark}

\newcommand{\mc}{\mathcal}

\usepackage{tikz}
\usetikzlibrary{calc}
\usetikzlibrary{positioning}
\usetikzlibrary{arrows.meta,shadows,positioning}
\usetikzlibrary{shapes.multipart}
\usetikzlibrary{backgrounds}
\tikzset{
  frame/.style={
    rectangle, draw,
    text width=6em, text centered,
    minimum height=4em,drop shadow,fill=white,
    rounded corners,
  },
  line/.style={
    draw, -{Latex},rounded corners=3mm,
  }
}

\newtheorem{remark}{Remark}[section]

\usepackage[english]{babel}

\author[a]{Kala Agbo Bidi}
\author[a]{Jean-Michel Coron}
\author[b]{Amaury Hayat}
\author[b,c]{Nathan Lichtlé}

\affil[a]{Laboratoire Jacques-Louis Lions, Sorbonne Universit\'{e}, Universit\'{e} de Paris, CNRS, INRIA, \'{e}quipe Cage, Paris, France (jean-michel.coron@sorbonne-universite.fr).}
\affil[b]{CERMICS, \'{E}cole des Ponts ParisTech, Champs-sur-Marne, France (amaury.hayat@enpc.fr, nathan.lichtle@enpc.fr).}
\affil[c]{Department of Electrical Engineering and Computer Science, UC Berkeley, Berkeley CA }

\date{\empty}

\title{Reinforcement Learning in Control Theory: A New Approach to Mathematical Problem Solving}
\begin{document}

\maketitle

\begin{abstract}
One of the central questions in control theory is achieving stability through feedback control. This paper introduces a novel approach that combines Reinforcement Learning (RL) with mathematical analysis to address this challenge, with a specific focus on the Sterile Insect Technique (SIT) system. The objective is to find a feedback control that stabilizes the mosquito population model. Despite the mathematical complexities and the absence of known solutions for this specific problem, our RL approach identifies a candidate solution for an explicit  stabilizing control.
This study underscores the synergy between AI and mathematics, opening new avenues for tackling intricate mathematical problems.

\end{abstract}

\section{Introduction}

AI for mathematics often refers to automatic theorem proving, either in formal language \cite{wu2020int,polu2022formal, HTPS, wu2022autoformalization}, or in natural language \cite{lewkowycz2022solving}, usually using language models. This article takes a different approach, presenting a Reinforcement Learning (RL) framework to solve a mathematical problem from control theory. The goal is to help mathematicians by finding a candidate solution to the problem 
in the spirit that checking a solution is often easier than finding it.

Control theory is about asking oneself: "if I can act on a system, what can I make it do?". In this area of mathematics, the system in question is usually a described by a set of differential (or partial differential) equations in which there is a component --called control-- that can be chosen. One of the main branches of this field, called stabilization, aims to find a way to make an equilibrium stable by choosing this control as a function of the state of the system. This is called a feedback control. Many mathematical techniques exist to solve this problem \cite{coron2007control}. However, in some cases the current mathematical theories are unable to find a successful feedback control. In this article we show that an RL approach combined with a mathematical analysis can help to find new mathematical feedback controls in such complicated cases.

We study a practical case: the SIT system, that models the control of insect pests (in particular mosquito population). SIT stands for \emph{Sterile Insect Technique} which consists in releasing sterilized insects to reduce or eliminate a target population. Initially used in agriculture to control  insect pests, it is today employed in the vector-born disease fight against mosquitoes
that carry illnesses such as  malaria and arboviruses \cite{almeida2019mosquito,alphey2011model} and there is a great interest both in research and in practice to understand which control to use \cite{bliman2018ensuring,bliman2019feedback,almeida2022optimal}. A more detailed overview of the literature is given in Appendix \ref{app:relatedworks}. From a mathematical point of view, without any control, this system of differential equations has a globally stable undesired equilibrium with high population of insect pests, and an unstable equilibrium with no insect pests, the \emph{zero equilibrium}. The mathematical goal is to find a feedback control such that the zero equilibrium is globally stable instead of the unwanted equilibrium.

The difficulties come from three reasons: showing global stability of dynamical systems (as opposed to local stability) is a very challenging mathematical problem for which there are only few mathematical tools; the system is not continuous, which is known in mathematics to bring some difficulties; we only have a partial measurement of the state of the system (see Section \ref{sec:framework} for more details). Because of these three difficulties, finding a feedback control for this system is a mathematically open question. With our approach we are able to derive an explicit mathematical feedback control to achieve the global stability. While there is still no mathematical proof that this feedback control is a solution to the problem, the numerical simulations 

strongly suggest that it is. 
We believe that this approach could be generalized to other open problems in control theory and give a new impulse to their resolution.

\section{Mathematical framework}

The dynamical system we consider is 
the SIT model for pests, here mosquitoes, given by
	\begin{align}
		&\dot{E} = \beta_E F \left(1-\frac{E}{K}\right) - \big( \nu_E + \delta_E \big) E,\label{eq:S1E1}  \\
		&\dot{M} = (1-\nu)\nu_E E - \delta_M M, \label{eq:S1E2} \\
		&\dot{F} =\nu\nu_E E \frac{M}{M+ M_s} - \delta_F F, \label{eq:S1E3} \\
		&\dot{M}_s = u - \delta_s M_s, \label{eq:S1E4}
	\end{align}
where
		 $E(t)\geq 0$ represents the mosquito density in aquatic phase,
		 $M(t)\geq 0$ the wild  adult  male density,
		 $F(t)\geq 0$ the density of adult feconded females,  $M_s(t)\geq 0$ the sterilized adult male density, and $u(t)\geq 0$, the control, is the density of sterilized
   males released at time $t$. We also denote the number of unfertilized 
   females by $F_{s}(t) =F(t)M_{s}(t)/M(t)$.

   When $u(t)= M_s(t) = 0$ for any $t\geq0$, the system \eqref{eq:S1E1}--\eqref{eq:S1E3} has a unique globally asymptotically stable equilibrium $(E(t),M(t),F(t)) \equiv (E^{*},M^{*},F^{*})$ where $E^{*}$, $M^{*}$ and $F^{*}$ are large constant values. This corresponds to the situation where mosquitoes reproduce freely. The state $(E(t),M(t),F(t)) \equiv(0,0,0)$ is also an equilibrium, albeit an unstable one. The mathematical problem is to find $u(t)$ of the form
   \begin{equation}
   \label{eq:control}
   u(t) = f(M(t)+M_{s}(t),F(t)+F_{s}(t)),
   \end{equation}
   where $f\in L^{\infty}(\mathbb{R}^{2})$ such that the zero equilibrium $(0,0,0)$ 
   is globally asymptotically stable and $M_s$ is asymptotically small, meaning
   there exists $c\in\mathbb{R}_+$ such that
      \begin{equation}  
      \label{eq:Ustar}
    \lim\limits_{t\rightarrow+\infty}\|u(t)\| = c< U^* := \frac{K\beta_E\nu(1-\nu)\nu_E^2\delta_s}{4(\delta_E+\nu_E)\delta_F\delta_M}\left(1-\frac{\delta_F(\nu_E+\delta_E)}{\beta_E\nu\nu_E}\right)^2,
   \end{equation}
and the equilibrium $(0,0,0,c/\delta_{s})$ of system \eqref{eq:S1E1}--\eqref{eq:S1E4} is globally asymptotically stable (see Definition \ref{def:GAS} below).
      
   \begin{defn}
   \label{def:GAS}
   The equilibrium $(0,0,0,c/\delta_{s})$ of the system \eqref{eq:S1E1}--\eqref{eq:S1E4} is \emph{globally asymptotically stable} if, for any initial condition $(E_{0}, M_{0},F_{0}, M_{s,0})$ there exists a unique solution $(E,M,F,M_{s})$ on $[0,+\infty)$ to the system \eqref{eq:S1E1}--\eqref{eq:S1E4} and for any $\varepsilon>0$ there exists $\delta>0$ such that
   \begin{gather}
\|(E_{0},M_{0},F_{0},M_{s,0}-c/\delta_{s})\|\leq \delta \implies \|(E(t),M(t),F(t), M_{s}(t)-c/\delta_{s})\|\leq \varepsilon,\;
\forall t\in[0,+\infty),\\
   \lim\limits_{t\rightarrow +\infty}\|(E(t),M(t),F(t),M_{s}(t)-c/\delta_{s})\| = 0,
   \end{gather}
   \end{defn}

The form constraint \eqref{eq:control} corresponds to a practical limitation: $M+M_{s}$ and $F+F_{s}$ are the total number of males and females which are typically what can be measured in practice 
(see \cite{almeida2019mosquito}).\\

\begin{remark}[Constant control]
The point of the constraint \eqref{eq:Ustar} is to avoid a constant control. Indeed, for a constant control $u(t) \equiv \bar{U}$, if $\bar{U} >U^{*}$ then the equilibrium $(0,0,0,U^{*}/\delta_{s})$ is globally asymptotically stable (see \cite{almeida2022optimal}).
In practice, one would like $c$ to be as small as possible in \eqref{eq:Ustar}. For the control we find in Section \ref{sec:results} with our approach, the value of $c$ is much smaller than $U^{*}$. In fact, in Appendix \ref{app:alternative}, we even show a simplified version of the control where $c$ can be chosen arbitrarily small.
\end{remark}

\section{Related works}

From a control theory perspective, several mathematical approaches have already been used in the literature 
to treat this problem 
either for the complete system \eqref{eq:S1E1}--\eqref{eq:S1E4} or for reduced models, 
using classical tools in control theory (control Lyapunov functions, LaSalle invariance principle, maximum principle, monotone dynamical systems, etc.
\cite{almeida2022optimal,2023-Agbo-Almeida-Coron-preprint,barclay1980sterile,anguelov2012mathematical,anguelov2020sustainable,2023Rossi}). See Appendix \ref{app:relatedworks} for more details. In particular, a feedback control was found in \cite{2023-Agbo-Almeida-Coron-preprint}, however this control depends on the four variables $(E,M,F,M_{s})$ and not only on the observable quantities $M+M_{s}$ and $F+F_{s}$.

Over the past few years, RL has emerged as a powerful approach for control, in a wide range of domains and applications (see Appendix \ref{app:relatedworks}). However, RL techniques, while powerful for decision-making, inherently provide control mechanisms that are discrete and numerical in nature. From a more rigorous mathematical point of view, these mechanisms often do not translate directly into analytical feedback control formulas. Aiming to address this limitation, in our work, we blend RL methodologies with mathematical analysis to extract an explicit mathematical control.

Using AI tools to help mathematicians by giving them an insight or a candidate solution was considered in \cite{davies2021advancing} using a different framework. Other approaches aimed to teach a model to guess mathematical solutions to a problem \cite{LampleCharton,DDSS}. However, in these approaches the solution of the mathematical problems involved are already known.

\section{Method}

\subsection{Our approach}

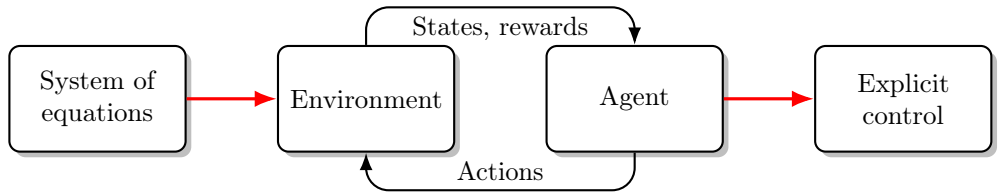
\begin{figure}
\centering
\begin{tikzpicture}[node distance=6cm, auto, background rectangle/.style={fill=white!15}, show background rectangle]
    % Block 1: System Equations
    \node [frame, thick] (model) {System of equations};

    % Block 2: RL Loop
    \node [frame, thick, right=1.2cm of model] (environment) {Environment};
    \node [frame, thick, right=1.2cm of environment] (agent) {Agent};

    \coordinate[below=1.2cm of environment] (P);
    \draw[thin,dashed] (P|-environment.west) -- (P|-environment.east);

    % Arrows for action, state, reward
    \draw[line,thick] (agent.south) -- ++ (0,-5mm) -| (environment.south) 
    node[above,pos=0.25,align=right] {Actions};
    \draw[line,thick] (environment.north) -- ++ (0,5mm) -| (agent.north) 
    node[below, pos=0.25, align=left] {States, rewards};

    \draw[line, very thick, color=red] (model.east) -- (environment.west) {};

    % Block 3: Explicit Controller Equations
    \node [frame, thick, right=1.2cm of agent] (explicit) {Explicit control};

    \draw[line, very thick, color=red] (agent.east) -- (explicit.west) {};
\end{tikzpicture}

\caption{Diagram representing the procedure by which we simulate our model in an environment that is used to train an RL agent, whose policy we then convert into an explicit control.}
\label{fig:diagram_method}
\end{figure}

Our proposed approach works in four steps, summarized in Figure \ref{fig:diagram_method}:

\begin{enumerate}
    \item[(Step 1)] \textbf{Discretize the equations} in a numerical scheme and use those dynamics to create a training environment by implementing the observations, actions, and rewards described in Section~\ref{sec:framework}.
    \item[(Step 2)] \textbf{Train an RL model} that learns to maximize the objective function we assign it through many simulations and obtain a numerical control feedback based on this numerical scheme.
    \item[(Step 3)] \textbf{Recover an explicit mathematical control} from the numerical control feedback. 
    \item[(Step 4)] \textbf{Perform several tests} using different numerical schemes and discretizations to ensure that the explicit control is efficient.
\end{enumerate}

\subsection{Reinforcement-Learning Framework}
\label{sec:framework}
Reinforcement Learning (RL) trains agents to optimize long-term rewards in various environments by maximizing the expected cumulative sum of rewards, denoted as $J(\pi_\theta)$, where $\pi_\theta$ is the policy, often guided by neural network weights $\theta$. Similar to minimizing a cost function $-J(\pi_\theta)$ in control theory, this can be modeled using a partially-observable Markov decision process (POMDP).
The observation space has two states: total males $M+M_s$ and females $F+F_s$, as the individual states $E, M, F$ and $M_s$ aren't independently measurable in the real world.
Observations are normalized to $[0,1]$. For better convergence, $M+M_s$ and $F+F_s$, which can typically range from $0$ to $100K$, are inputted to the neural network at varied scales and normalized. The single action $a_t$ ranges from [-1,1], remapped to $u(t) \in [0, 10K]$ for model equations \eqref{eq:S1E1}--\eqref{eq:S1E4}. To simplify training by artificially reducing the horizon, each action is repeated for multiple simulation steps.

Finally, our optimization criterion takes the following form at time step $t$:
$$r_t = c_1 \lVert E(t), M(t), F(t) \rVert_2 + c_2(t) \lVert M_s(t) \rVert_2 \quad \quad \text{ with } \quad \quad c_2(t) = \begin{cases}
    c_3 &\text{ if } t < 0.9T, \\
    c_3 + c_4 &\text{ otherwise.}
\end{cases} $$  
Near the horizon's end, specifically when $0.9T \leq t \leq T$, we introduce a positive weight $c_4$ to the penalty on $\lVert M_s(t) \rVert_2$. This aids steer the RL training toward the desired asymptotic convergence of the state.
Parameters and more experiment details can be found in Appendix \ref{app:expdetails}.

\section{Main Results}
\label{sec:results}
The trained RL model converges after around 10M steps to a numerical control that we represent in Figure~\ref{fig:nn_heatmap_mms_ffs} (see Appendix \ref{app:heatmap}) as a function of $M+M_{s}$ and $F+F_{s}$. We see that the plot of the control in linear scale is not really informative (see Fig. \ref{fig:regression} left). However, in log scales the expression of the control seems clearer (see Fig. \ref{fig:regression} right) and clearly has two parts. In each of them the control seems to be close to a bang-bang control with a thin transition. With a simple regression we approximate this numerical control with the explicit mathematical control
\begin{equation}
    u_\text{reg}(M+M_s, F+F_s) = \begin{cases}
        u_\text{reg}^\text{left}(M+M_s, F+F_s) & \text{ if } M+M_s < M^{*}, \\
        u_\text{reg}^\text{right}(M+M_s, F+F_s) & \text{ otherwise,}
    \end{cases}
    \label{eq:controlv1}
\end{equation}
where $u$ is defined on $(0,+\infty)^{2}$ and 

\begin{minipage}[t]{0.45\textwidth}
{\footnotesize
\begin{align*}
    u_\text{reg}^\text{left} &= \begin{cases}u_\text{min} & \text{ if } I_{1}(F+F_{s}) > \alpha_{2}, \\ u_\text{max} \left(\alpha_{2}-I_{1}
    \right) & \text{ if } I_{1}\in(\alpha_{1},\alpha_{2}], \\ u_\text{max} & \text{ otherwise, and} \end{cases} \\
    \end{align*}}
\end{minipage}
\begin{minipage}[t]{0.45\textwidth}
{\footnotesize    \begin{align*}
    u_\text{reg}^\text{right} &= \begin{cases}
u_\text{min} & \text{ if } I_{2}>\alpha_{2}, \\ u_\text{max} \left(  \alpha_{2} - I_{2} \right) & \text{ if } I_{2}\in(\alpha_{1},\alpha_{2}], \\ u_\text{max} & \text{ otherwise.}
\end{cases}
\end{align*}
}
\end{minipage}

where $I_{1}(x) = \frac{\log (M^{*})}{\log(x)}$ and $I_{2}(x,y) =  \frac{\log(x)}{\log(y)}$,  $M^{*}=200$, $\alpha_{1} = 3$, $\alpha_{2} = 4$, $u_{\max}=3 \cdot 10^{5}$ is imposed by physical constraints and $u_{\min}$ can be chosen.

During training and evaluation the numerical feedback control includes a slight noise to enhance robustness and exploration. Surprisingly, when tested with $u_{\min}=0$, the mathematical control \eqref{eq:controlv1} with a small additional noise $\eta(t)$ exhibits asymptotic stability, whereas noise-free control does not, displaying cyclic-like behavior (see Figure \ref{fig:regression}). The paradox arises from $u_{\min}=0$. Introducing a small positive value $\varepsilon$ for this parameter is enough to obtain asymptotic stability. In the noisy control, because of the condition $u\geq 0$, $u_{\min}+\eta(t)$ is positive in average, explaining the stabilization. More details are given in Appendix \ref{app:noisy}.

The efficiency of the feedback control \eqref{eq:controlv1} is illustrated in Appendix \ref{app:results} on many numerical simulations for a large array of initial conditions with different 
discretizations,
suggesting that this nonlinear control is indeed a solution to the mathematical problem considered. With this control, the system converges globally to the equilibrium $(E^{*},M^{*},F^{*},M_{s}^{*}) = (0,0,0,\lambda)$ where $\lambda = \varepsilon/\delta_{s}$ and $\varepsilon$ can be chosen much smaller than the $U^{*}$ given in \eqref{eq:Ustar}. 
A more detailed analysis of the result is given in Appendix \ref{app:results}.

\section{Discussion and conclusion}
We presented an RL framework to solve a type of mathematical problem in control theory, and we used it to find an explicit candidate solution for the stabilization of the SIT system. In the future, it would be interesting to use a regression to explicit the reward $J$ as a function of the initial state in order to obtain an explicit Lyapunov function which would show the asymptotic stability of the system with the explicit control.
Going further, this approach could be likely generalized to other systems. This is an incentive to use more AI techniques to solve mathematical problems, especially in control theory.

\bibliographystyle{plain}
\bibliography{main_arxiv}

\begin{thebibliography}{10}

\bibitem{2023-Agbo-Almeida-Coron-preprint}
Kala Agbo~Bidi, Luis Almeida, and Jean-Michel Coron.
\newblock Global stabilization of sterile insect technique model by feedback
  laws.
\newblock {\em arXiv, 2307.00846}, 2023.

\bibitem{almeida2019mosquito}
Luis Almeida, Michel Duprez, Yannick Privat, and Nicolas Vauchelet.
\newblock Mosquito population control strategies for fighting against
  arboviruses.
\newblock {\em Mathematical Biosciences and Engineering}, 16(6):6274--6297,
  2019.

\bibitem{almeida2022optimal}
Lu{\'\i}s Almeida, Michel Duprez, Yannick Privat, and Nicolas Vauchelet.
\newblock Optimal control strategies for the sterile mosquitoes technique.
\newblock {\em Journal of Differential Equations}, 311:229--266, 2022.

\bibitem{alphey2011model}
Nina Alphey, Luke Alphey, and Michael~B Bonsall.
\newblock A model framework to estimate impact and cost of genetics-based
  sterile insect methods for dengue vector control.
\newblock {\em PLoS One}, 6(10):e25384, 2011.

\bibitem{anguelov2012mathematical}
Roumen Anguelov, Yves Dumont, and Jean Lubuma.
\newblock Mathematical modeling of sterile insect technology for control of
  anopheles mosquito.
\newblock {\em Computers \& Mathematics with Applications}, 64(3):374--389,
  2012.

\bibitem{anguelov2020sustainable}
Roumen Anguelov, Yves Dumont, and Ivric~Valaire Yatat~Djeumen.
\newblock Sustainable vector/pest control using the permanent sterile insect
  technique.
\newblock {\em Mathematical Methods in the Applied Sciences},
  43(18):10391--10412, 2020.

\bibitem{barclay1980sterile}
H~Barclay and M~Mackauer.
\newblock The sterile insect release method for pest control: a
  density-dependent model.
\newblock {\em Environmental Entomology}, 9(6):810--817, 1980.

\bibitem{bliman2019feedback}
Pierre-Alexandre Bliman.
\newblock Feedback control principles for biological control of dengue vectors.
\newblock In {\em 2019 18th European Control Conference (ECC)}, pages
  1659--1664. IEEE, 2019.

\bibitem{bliman2018ensuring}
Pierre-Alexandre Bliman, M~Soledad Aronna, Fl{\'a}vio~C Coelho, and Moacyr~AHB
  da~Silva.
\newblock Ensuring successful introduction of wolbachia in natural populations
  of aedes aegypti by means of feedback control.
\newblock {\em Journal of mathematical biology}, 76:1269--1300, 2018.

\bibitem{2019-Bliman-Cardona-Salgado-Dumont-Vasilieva-MB}
Pierre-Alexandre Bliman, Daiver Cardona-Salgado, Yves Dumont, and Olga
  Vasilieva.
\newblock Implementation of control strategies for sterile insect techniques.
\newblock {\em Math. Biosci.}, 314:43--60, 2019.

\bibitem{2022-Bliman-Dumont-MB}
Pierre-Alexandre Bliman and Yves Dumont.
\newblock Robust control strategy by the {S}terile {I}nsect {T}echnique for
  reducing epidemiological risk in presence of vector migration.
\newblock {\em Math. Biosci.}, 350:Paper No. 108856, 23, 2022.

\bibitem{DDSS}
Francois Charton, Amaury Hayat, and Guillaume Lample.
\newblock Learning advanced mathematical computations from examples.
\newblock In {\em International Conference on Learning Representations}, 2020.

\bibitem{coron2007control}
Jean-Michel Coron.
\newblock {\em Control and nonlinearity}.
\newblock American Mathematical Soc., 2007.

\bibitem{2023Rossi}
Andrea Cristofaro and Luca Rossi.
\newblock {Backstepping control for the sterile mosquitoes technique:
  stabilization of extinction equilibrium}.
\newblock working paper or preprint, 2023.

\bibitem{2023-Cristofaro-Rossi}
Andrea Cristofaro and Luca Rossi.
\newblock {Backstepping control for the sterile mosquitoes technique:
  stabilization of extinction equilibrium}.
\newblock {\em Preprint}, 2023.

\bibitem{davies2021advancing}
Alex Davies, Petar Veli{\v{c}}kovi{\'c}, Lars Buesing, Sam Blackwell, Daniel
  Zheng, Nenad Toma{\v{s}}ev, Richard Tanburn, Peter Battaglia, Charles
  Blundell, Andr{\'a}s Juh{\'a}sz, et~al.
\newblock Advancing mathematics by guiding human intuition with ai.
\newblock {\em Nature}, 600(7887):70--74, 2021.

\bibitem{7963427}
Amir-massoud Farahmand, Saleh Nabi, and Daniel~N. Nikovski.
\newblock Deep reinforcement learning for partial differential equation
  control.
\newblock In {\em 2017 American Control Conference (ACC)}, pages 3120--3127,
  2017.

\bibitem{7989385}
Shixiang Gu, Ethan Holly, Timothy Lillicrap, and Sergey Levine.
\newblock Deep reinforcement learning for robotic manipulation with
  asynchronous off-policy updates.
\newblock In {\em 2017 IEEE International Conference on Robotics and Automation
  (ICRA)}, pages 3389--3396, 2017.

\bibitem{kiumarsi2017optimal}
Bahare Kiumarsi, Kyriakos~G Vamvoudakis, Hamidreza Modares, and Frank~L Lewis.
\newblock Optimal and autonomous control using reinforcement learning: A
  survey.
\newblock {\em IEEE transactions on neural networks and learning systems},
  29(6):2042--2062, 2017.

\bibitem{LampleCharton}
Guillaume Lample and Fran{\c{c}}ois Charton.
\newblock Deep learning for symbolic mathematics.
\newblock In {\em International Conference on Learning Representations}, 2019.

\bibitem{HTPS}
Guillaume Lample, Marie-Anne Lachaux, Thibaut Lavril, Xavier Martinet, Amaury
  Hayat, Gabriel Ebner, Aur{\'e}lien Rodriguez, and Timoth{\'e}e Lacroix.
\newblock {HyperTree Proof Search for Neural Theorem Proving}.
\newblock {\em {A}dvances in neural information processing systems}, 2022.

\bibitem{lewkowycz2022solving}
Aitor Lewkowycz, Anders Andreassen, David Dohan, Ethan Dyer, Henryk
  Michalewski, Vinay Ramasesh, Ambrose Slone, Cem Anil, Imanol Schlag, Theo
  Gutman-Solo, et~al.
\newblock Solving quantitative reasoning problems with language models.
\newblock {\em Advances in Neural Information Processing Systems},
  35:3843--3857, 2022.

\bibitem{lillicrap2015continuous}
Timothy~P Lillicrap, Jonathan~J Hunt, Alexander Pritzel, Nicolas Heess, Tom
  Erez, Yuval Tassa, David Silver, and Daan Wierstra.
\newblock Continuous control with deep reinforcement learning.
\newblock {\em arXiv preprint arXiv:1509.02971}, 2015.

\bibitem{mahmood2018benchmarking}
A.~Rupam Mahmood, Dmytro Korenkevych, Gautham Vasan, William Ma, and James
  Bergstra.
\newblock Benchmarking reinforcement learning algorithms on real-world robots,
  2018.

\bibitem{mnih2013playing}
Volodymyr Mnih, Koray Kavukcuoglu, David Silver, Alex Graves, Ioannis
  Antonoglou, Daan Wierstra, and Martin Riedmiller.
\newblock Playing atari with deep reinforcement learning.
\newblock {\em arXiv preprint arXiv:1312.5602}, 2013.

\bibitem{polu2022formal}
Stanislas Polu, Jesse~Michael Han, Kunhao Zheng, Mantas Baksys, Igor
  Babuschkin, and Ilya Sutskever.
\newblock Formal mathematics statement curriculum learning.
\newblock In {\em 11th International Conference on Learning Representations},
  2022.

\bibitem{raffin2019stable}
Antonin Raffin, Ashley Hill, Adam Gleave, Anssi Kanervisto, Maximilian
  Ernestus, and Noah Dormann.
\newblock Stable-baselines3: Reliable reinforcement learning implementations.
\newblock {\em Journal of Machine Learning Research}, 22(268):1--8, 2021.

\bibitem{schulman2017ppo}
John Schulman, Filip Wolski, Prafulla Dhariwal, Alec Radford, and Oleg Klimov.
\newblock Proximal policy optimization algorithms.
\newblock {\em arXiv preprint arXiv:1707.06347}, 2017.

\bibitem{silver2018general}
David Silver, Thomas Hubert, Julian Schrittwieser, Ioannis Antonoglou, Matthew
  Lai, Arthur Guez, Marc Lanctot, Laurent Sifre, Dharshan Kumaran, Thore
  Graepel, et~al.
\newblock A general reinforcement learning algorithm that masters chess, shogi,
  and go through self-play.
\newblock {\em Science}, 362(6419):1140--1144, 2018.

\bibitem{strugarek2019use}
Martin Strugarek, Herv{\'e} Bossin, and Yves Dumont.
\newblock On the use of the sterile insect release technique to reduce or
  eliminate mosquito populations.
\newblock {\em Applied Mathematical Modelling}, 68:443--470, 2019.

\bibitem{vinyals2019grandmaster}
Oriol Vinyals, Igor Babuschkin, Wojciech~M Czarnecki, Micha{\"e}l Mathieu,
  Andrew Dudzik, Junyoung Chung, David~H Choi, Richard Powell, Timo Ewalds,
  Petko Georgiev, et~al.
\newblock Grandmaster level in starcraft ii using multi-agent reinforcement
  learning.
\newblock {\em Nature}, 575(7782):350--354, 2019.

\bibitem{wu2020int}
Yuhuai Wu, Albert~Qiaochu Jiang, Jimmy Ba, and Roger Grosse.
\newblock Int: An inequality benchmark for evaluating generalization in theorem
  proving.
\newblock {\em arXiv preprint arXiv:2007.02924}, 2020.

\bibitem{wu2022autoformalization}
Yuhuai Wu, Albert~Qiaochu Jiang, Wenda Li, Markus Rabe, Charles Staats, Mateja
  Jamnik, and Christian Szegedy.
\newblock Autoformalization with large language models.
\newblock {\em {Advances in Neural Information Processing Systems}},
  35:32353--32368, 2022.

\end{thebibliography}

\newpage
\appendix
\section{Mathematical system interpretation and parameters}

   In system \eqref{eq:S1E1}--\eqref{eq:S1E4}, we assumed that all females are immediately fertilized when they emerge from the pupal stage. The equation on $F$ makes sense  when we add the sterile male in which case only a fraction of the females will be fertilized. The interpretation of the parameters are given below \cite{almeida2022optimal}:
   \begin{itemize}
    \item $\beta_E>0$ is the oviposition rate,
    \item $\delta_E,\delta_M,\delta_F >0 $ are the death rates for eggs, wild adult males and fertilized females respectively,
    \item $\nu_E>0$ is the hatching rate for eggs,
    \item $\nu\in (0,1)$ the probability that a pupa gives rise to a female (and $(1-\nu)$ is, therefore, the probability to give rise to a male),
    \item $\delta_s>0$ is the death rate of sterilized adult, 
    \item $K>0$ is the environmental capacity for eggs. It can be interpreted as the maximum density of eggs that females can lay in breeding sites. Since here the larval and pupal compartments are not present, it is as if $E$ represents all the aquatic compartments in which case in this term $K$ represents a logistic law's carrying capacity for the aquatic phase that also includes the effects of competition between larvae.
\end{itemize}

Besides, we also assume that $\delta_s\geq \delta_M$, which is usually considered as a biologically relevant assumption \cite{almeida2022optimal} 
Typical values for these parameters can be found in \cite{strugarek2019use} and are given in Table \ref{tab:tableparams}.

\begin{table}[ht]
    \centering
    \begin{tabular}{|c|c|c|c|c|}
        \hline
        Parameter & Name & Value Interval & Chosen Value & Unity \\
        \hline
        $\beta_E$ & Effective fecundity & [7.46, 14.85] & 8 & Day$^{-1}$ \\
        \hline
        $\nu_E$ & Hatching parameter & [0.005, 0.25] & 0.25 & Day$^{-1}$ \\
        \hline
        $\delta_E$ & Aquatic phase death rate & [0.023, 0.046] & 0.03 & Day$^{-1}$ \\
        \hline
        $\delta_F$ & Female death rate & [0.033, 0.046] & 0.04 & Day$^{-1}$ \\
        \hline
        $\delta_M$ & Males death rate & [0.077, 0.139] & 0.1 & Day$^{-1}$ \\
        \hline
        $\delta_s$ & Sterilized male death rate & - & 0.12 & Day$^{-1}$ \\
        \hline
        $\nu$ & Probability of emergence & - & 0.49 & \\
        \hline
        K & Environmental capacity for eggs & - & 50000 & \\
        \hline
    \end{tabular}
    \vspace{0.2cm}
    \caption{Parameters for the system \eqref{eq:S1E1}--\eqref{eq:S1E4}.}
    \label{tab:tableparams}
\end{table}

\section{Experiment details}
\label{app:expdetails}
We train our RL policies using proximal policy optimization (PPO) \cite{schulman2017ppo}, a state-of-the-art policy gradient algorithm. We use the implementation of PPO provided in Stable Baselines 3 \cite{raffin2019stable} (version 1.6.2, Python 3.8), a popular RL library that provides a collection of state-of-the-art algorithm implementations, as well as various tools for RL research.

The models are trained for 10 million environment timesteps (or 7 billion simulation timesteps) on 12 CPUs, which takes about 7 hours. During each iteration, we collect 12288 (1024 per CPU) environment steps, then run 5 epochs of optimization with a batch size of 1024. The agent's policy is a fully-connected neural network with 2 hidden layers of 256 neurons each, with $\tanh$ non-linearities between each layer, outputting the mean and standard deviation of a normal distribution that is then used to sample the action. More formally, for a given observation vector, the neural network policy outputs a mean $\mu_t$ and standard deviation $\sigma_t$ and the action is sampled as $a_t \sim \mc N(\mu_t, \sigma_t)$. We train with a learning rate of $3 \times 10^{-4}$, gamma factor $\gamma=0.99$, and all other hyperparameters are left to their default values.

We run each simulation for $T=1001$ days (or 143 weeks), with a timestep $dt = 0.01$ days, and each action is repeated $700$ times, meaning that the environment horizon is $143$ steps and a new action is taken each week. For each simulation, the initial condition is uniformly sampled between $0$ and $10K$: $E(0), M(0), F(0), M_s(0) \sim \mc U(0, 10K)$. For our reward function, we use coefficients $c_1 = 0.1$, $c_3 = 0.001$ and $c_4 = 0.01$. 

\section{Related works}
\label{app:relatedworks}

From a mathematical point of view, several mathematical techniques have been used, either for this model or reduced models. In particular, two reduced model have been considered: A two dimensional model (2D-model) obtained by assuming that the dynamics of males and eggs are fast so that these two populations can be assumed to be at equilibrium (see \cite[$(\mathcal{S}_1)$, page 231-232]{almeida2022optimal} or \cite[(2)]{2023-Cristofaro-Rossi}) and a three dimensional model (3D-model) obtained by  overlooking  the non-adult stages (see\cite[(7a)-(7b)-(7c)]{2019-Bliman-Cardona-Salgado-Dumont-Vasilieva-MB}). 
These mathematical approaches have led to the following stabilizing feedback controls:
\begin{itemize}
\item Stabilization using impulsive feedback controls for the 3D-model:   \cite[Theorem  6]{2019-Bliman-Cardona-Salgado-Dumont-Vasilieva-MB}, and \cite[Theorem  7]{2019-Bliman-Cardona-Salgado-Dumont-Vasilieva-MB} for the case of  sparse measurements. The case of vector migration is also considered in \cite{2022-Bliman-Dumont-MB}.
\item Stabilization using optimal feedback controls for the 2D-model: \cite[Remark 4]{almeida2022optimal}.
\item Stabilization using the backstepping method: see \cite{2023-Cristofaro-Rossi} for the 2D-model and \cite[Section 3.1]{2023-Agbo-Almeida-Coron-preprint} for the complete model.
\item Stabilization using simple linear feedback laws for which the stabilization for the complete model is conjectured and proved for positively invariant subsets \cite[Sections 3.2 and 3.3]{2023-Agbo-Almeida-Coron-preprint}.
\end{itemize}

Our approach differs by using deep reinforcement learning to construct control feedback laws. In the past few years, RL has emerged as a powerful approach for control, leveraging its ability to learn near-optimal decision making strategies through interactions with an environment, and has demonstrated remarkable successes across a wide range of domains and applications. 
In robotics, RL has enabled machines to learn complex control tasks such as locomotion, manipulation, and dexterous object handling \cite{7989385, mahmood2018benchmarking,lillicrap2015continuous}. In the realm of games, RL algorithms have achieved superhuman performance in challenging domains like Go, chess, and StarCraft \cite{silver2018general, vinyals2019grandmaster}. Moreover, RL has excelled in playing classic Atari games, surpassing human-level performance by learning directly from pixel inputs \cite{mnih2013playing}. 
RL can also excel in controlling ODE/PDE problems directly: \cite{kiumarsi2017optimal} applies RL to optimal control problems, \cite{7963427} applies RL to a flow control problem modeled by PDEs.
However our approach differs from these by fitting an explicit mathematical control law using the learned neural network. 
These remarkable achievements highlight the versatility and potential of RL as a general-purpose approach for solving complex practical control problems in diverse domains.

\section{Results}
\label{app:results}
\subsection{Numerical control after RL training}
\label{app:heatmap}
\begin{figure*}[h!]
\centering
\begin{subfigure}{.5\textwidth}
  \centering
  \includegraphics[width=.99\linewidth]{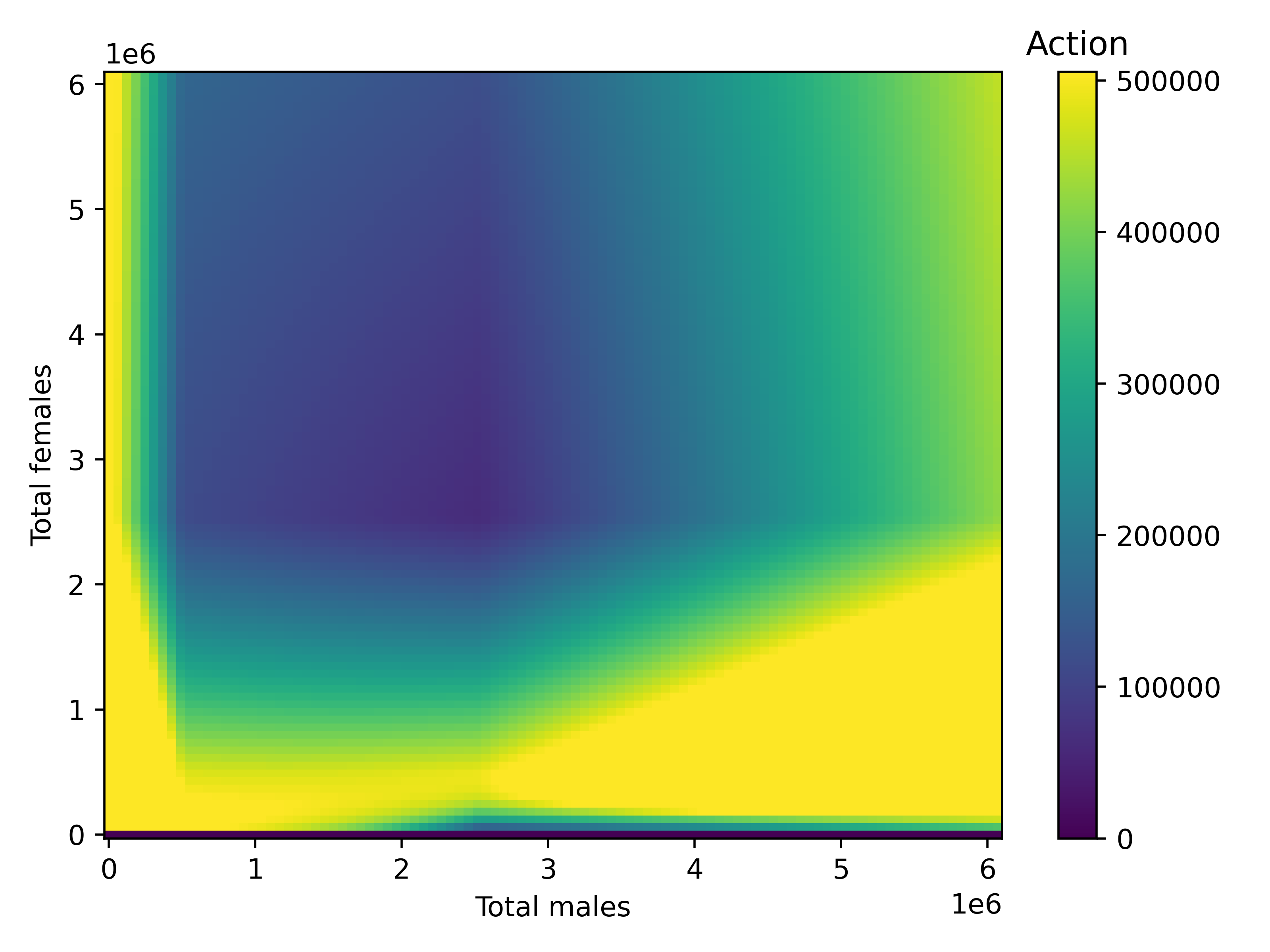}
\end{subfigure}%
\begin{subfigure}{.5\textwidth}
  \centering
  \includegraphics[width=.99\linewidth]{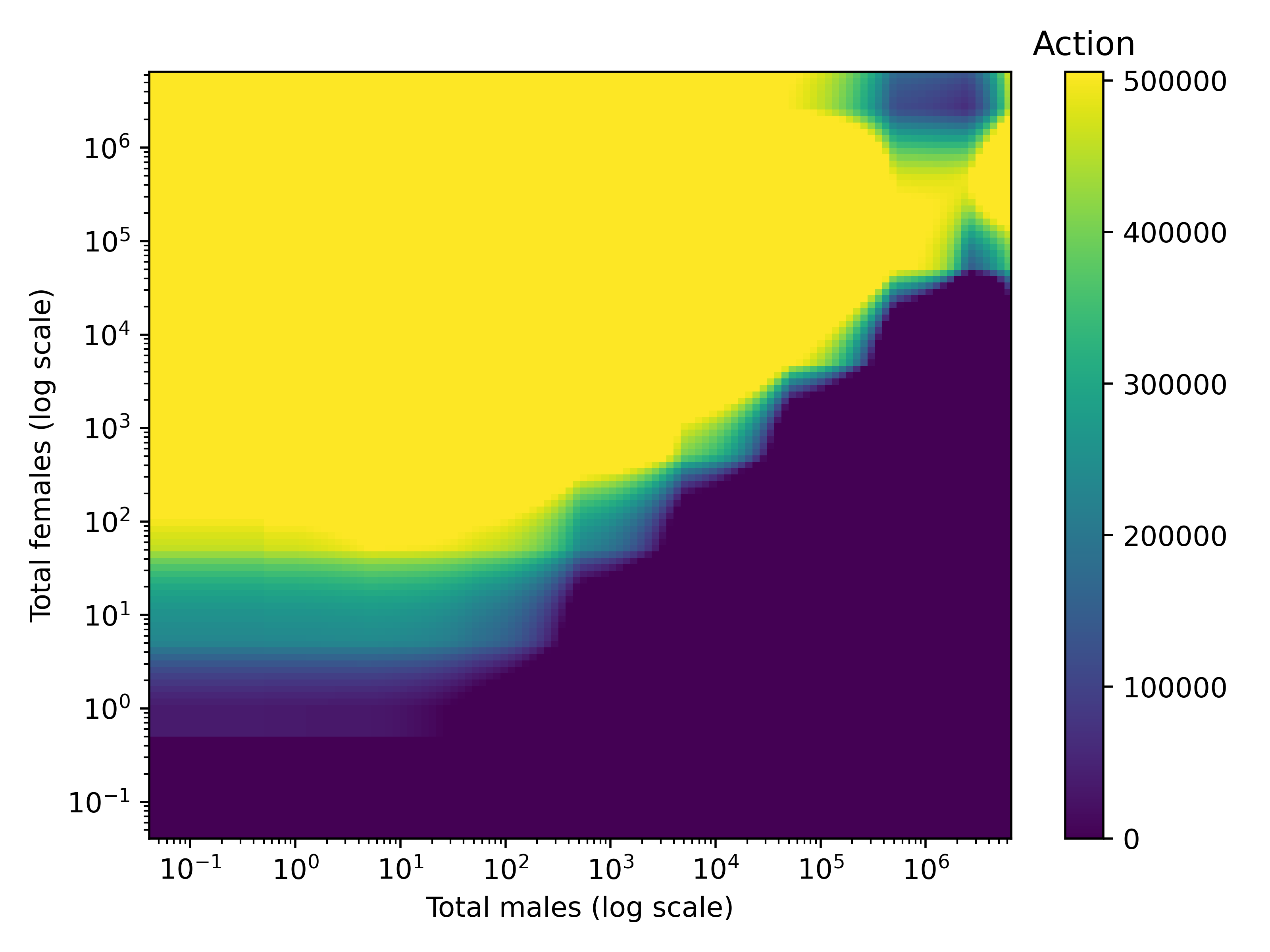}
\end{subfigure}
\caption{Heatmap of the model's action $u$ as a function of $M+M_s$ and $F+F_s$, in linear scale (left) and logarithmic scale (right).
\label{fig:nn_heatmap_mms_ffs}}

\end{figure*}

In Figure \ref{fig:nn_heatmap_mms_ffs} we represent the numerical control, that is the model's action as a function of $M+M_{s}$ (total males) and $F+F_{s}$ (total females).
Interestingly, the plot in linear scale (left) is not very informative, suggesting that it would be complicated to have a good regression directly as a function of $M+M_s$ and $F+F_s$. However, in log scale (right) the form of the function seems much more identified.

\subsection{Effect of the noise}
\label{app:noisy}

\begin{figure*}[h!]
\centering
\begin{subfigure}{.5\textwidth}
  \centering
  \includegraphics[width=.99\linewidth]{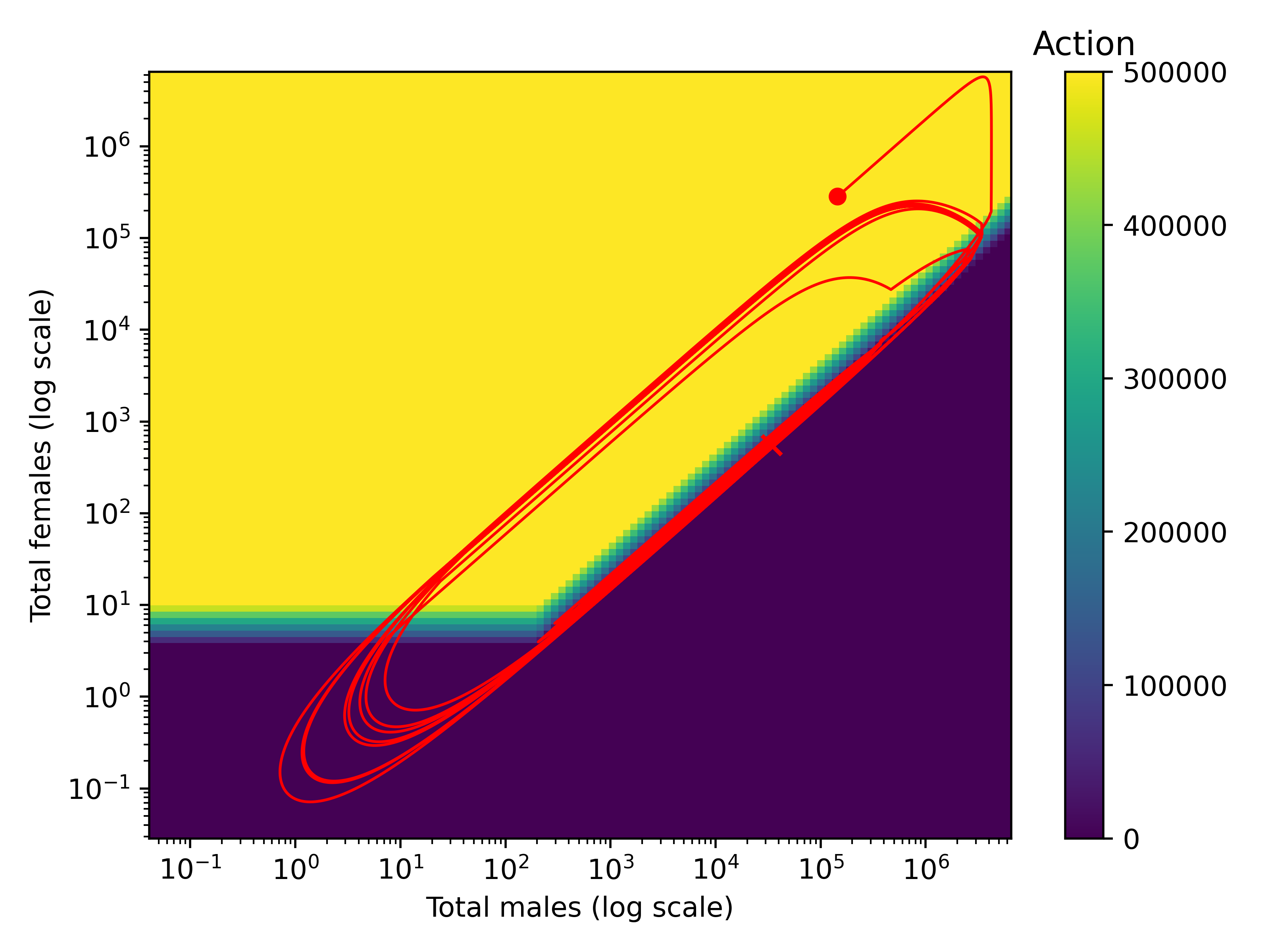}
\end{subfigure}%
\begin{subfigure}{.5\textwidth}
  \centering
  \includegraphics[width=.99\linewidth]{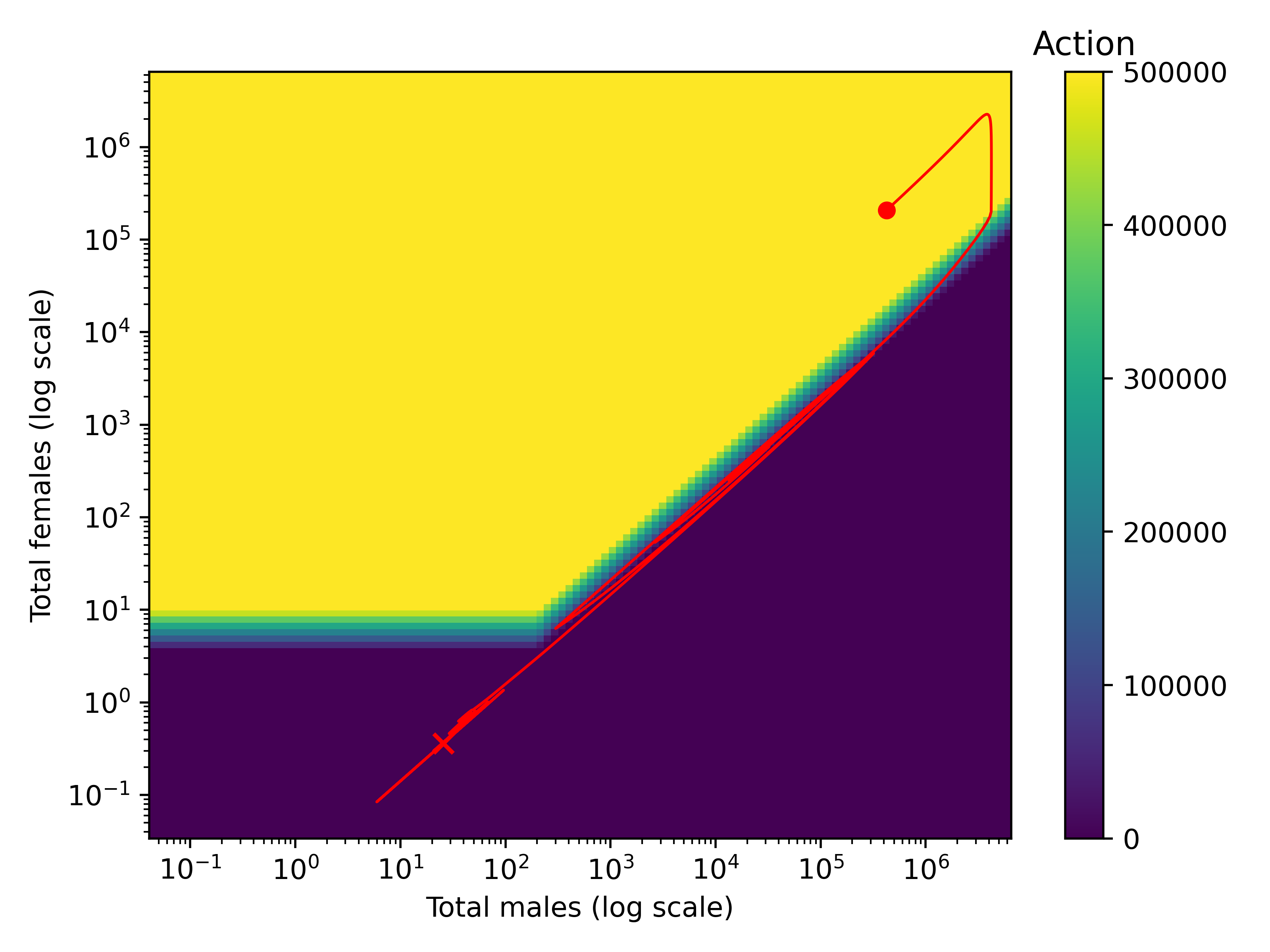}
\end{subfigure}
\caption{Heatmap of the regression model's action $u(M+M_s, F+F_s)$ as a function of total males and total females. A state-space trajectory is plotted in red, with the dot indicating initial state and the cross final state, for the heatmap only (left) and when a small noise $\mu \sim \mathcal N(0, 5)$ is added on top of the action (right).
\label{fig:regression}}
\end{figure*}

During training and testing the numerical control feedback law includes by default a small noise. This ensures some robustness of the control and a good exploration. We tested the mathematical control we derived (given in \eqref{eq:controlv1}) with and without noise. To our surprise, the control with a small noise does seem to ensure the asymptotic stability, whereas the control without any noise does not seem to. Indeed, without noise, the control seems to have a cyclic behavior and never converges (see Figure \ref{fig:regression} (left)). When adding a small noise, however, the stability is restored (see Figure \ref{fig:regression} (right)). The explication to this 
apparent paradox is that having exactly $u_{\min}=0$ in the one of the branch of the control given in \eqref{eq:controlv1} is apparently too strong to allow the model to converge completely. Replacing the value with $u_{\min}>0$ for a small $u_{\min}$ (typically $u_{\min}=10$) allows to stabilize the system without noise (see Figure \ref{fig:regression} (right)). In the system with noise, because the control $u$ has to be positive, the noise increases in average the effective value of $u_{\min}$ of the control \eqref{eq:controlv1}, which explains the apparent stabilization. Note that 
 $(0,0,0,u_{\min}/\delta_s)$ is an equilibrium of the system stabilized, which solves the problem described in Section \ref{sec:results} provided that $u_{\min}<U^{*}$. Here, with the values of Table \ref{tab:tableparams}, 
  $$U^{*}\approx 1.6\;10^{5}\text{ while }u_{\min}=10,$$ 
  which means that this control is a very good solution to the problem.

\subsection{Effectiveness of the mathematical control}

\begin{figure}
    \centering
    \includegraphics[width=0.99\linewidth]{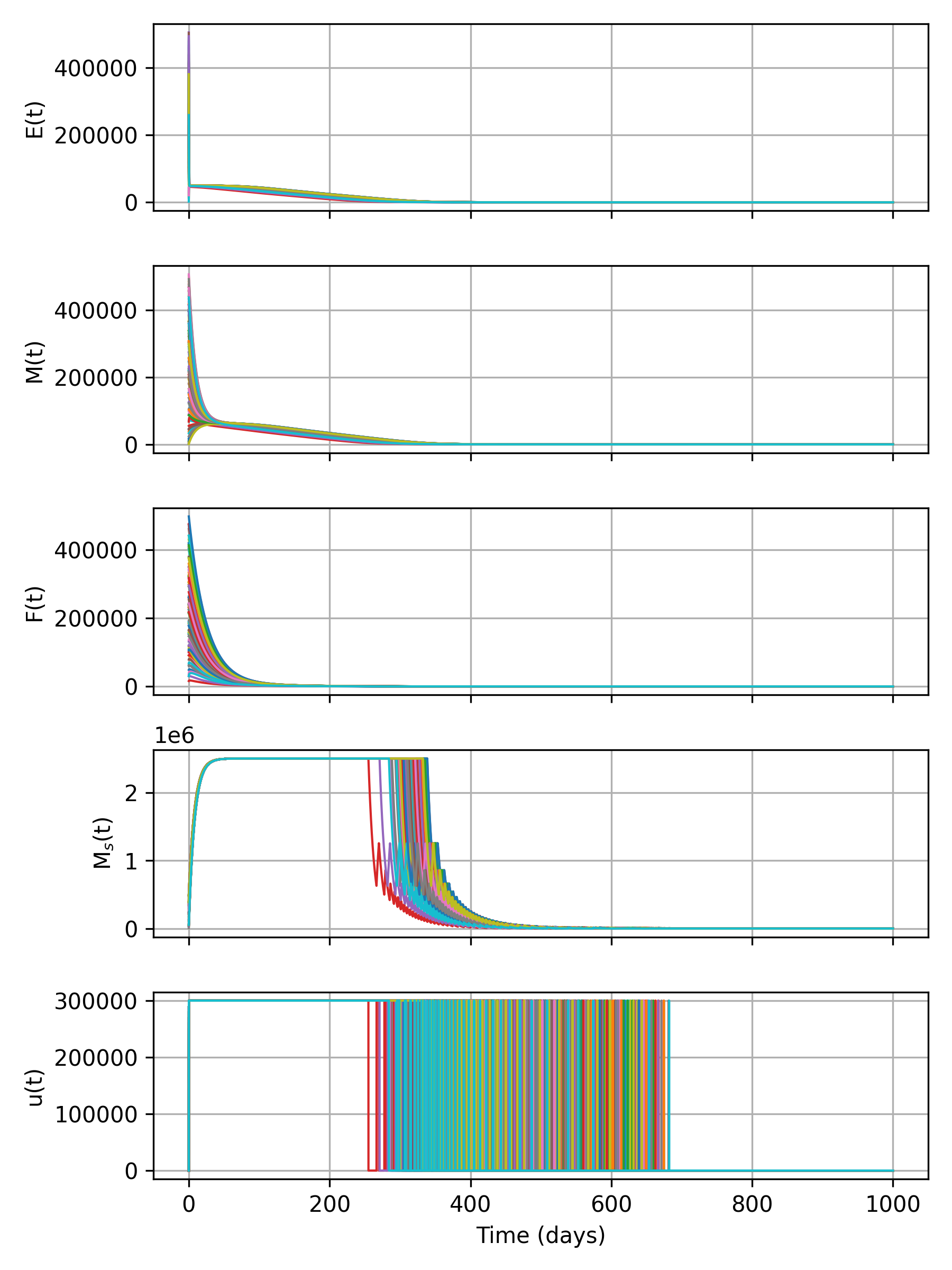}
    \caption{States and control $v_\text{reg}$ with $u_\text{min}=5$ and $u_\text{max}=300000$ over a duration of 1000 days for 100 simulations with random initial conditions in $[0, 10K]^4$. Each color correspond to a simulation.}
    \label{fig:100sims_reg2}
\end{figure}

We test the candidate mathematical control \eqref{eq:controlv1} with random initial condition $(E_{0},M_{0}, F_{0}, M_{s,0})$ in $[0,10K]$ 
to check the stability of the equilibrium $(0,0,0,u_{\min}/\delta_{s})$. This is represented in Figure \ref{fig:100sims_reg2}.
 We also show the numerical values obtained and their variance in Table \ref{fig:tableMMS}.

\begin{table}
    \centering
    \begin{tabular}{|l |c| c| c| c|}
        \hline
         & 200 days & 400 days & 600 days & 800 days \\
        \hline
        average $|E|+|M|+|F|$ & 50,801.64 & 8,020.47 & 96.25 & 0.60 \\
        \hline
        variance $|E|+|M|+|F|$ & 45,422,159 & 4,012,119 & 1,394 & 0.02 \\
        \hline
        maximum $|E|+|M|+|F|$ & 59,026.78 & 10,455.31 & 149.03 & 0.80 \\
        \hline
        average $|M_{s}|$ & 2,207,795.88 & 693,675.84 & 16,743.12 & 50.27 \\
        \hline
        variance $|M_{s}|$ & 53,101,242,728 & 13,102,504,436 & 41,355,608 & 81.24 \\
        \hline
        maximum $|M_{s}|$ & 2,473,954.23 & 822,154.59 & 25,783.57 & 70.59 \\
        \hline
    \end{tabular}
    \caption{Statistics over 100 simulations with random initial conditions in $[0, 10K]^4$ using control $u_\text{reg}$ (see \eqref{eq:controlv1}) with $u_\text{min}=5$ and $u_\text{max}=300000$ over a duration of 600 days.}
    \label{fig:tableMMS}
\end{table}

\section{An alternative control}
\label{app:alternative}
\begin{figure}
    \centering
    \includegraphics[width=0.99\textwidth]{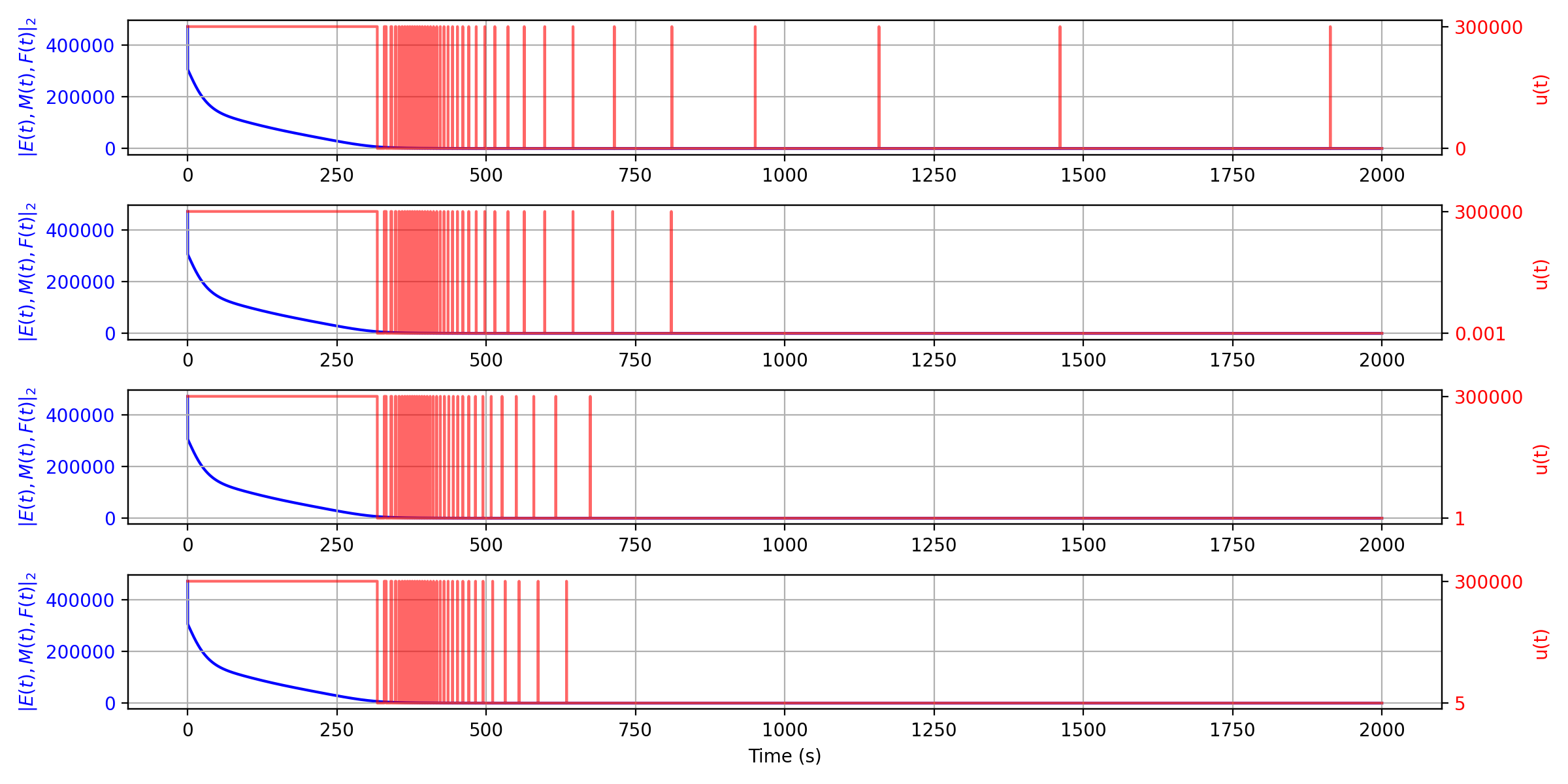}
    \caption{Norm of the states $\lVert E(t), M(t), F(t) \rVert_2$ (blue) and control $u(t)$ (red) as a function of time for different values of $u_\text{min}$ (0, 0.001, 1, and 5 respectively from top to bottom) and $u_\text{max} = 300,000$, over 2000 days and with the same initial condition.}
    \label{fig:controlv2_different_umin}
\end{figure}

Using the candidate solution provided by the method, we simplified \eqref{eq:controlv1} and obtained a second candidate feedback control given by
\begin{equation}
v_\text{reg}(M+M_s, F+F_s) = \begin{cases}
        u_\text{min} & \text{ if } \frac{\log (M+M_s)}{\log (F+F_s)} > \alpha_{2}, \\
        u_\text{max} & \text{ otherwise,}
    \end{cases}
\label{eq:controlv2}
\end{equation}
where $\alpha_{2}=4$ with the parameters of Table \ref{tab:tableparams} and $u_{\max}=300000$ is still imposed.
Interestingly, not only is this control faster to converge (compare Tables \ref{fig:tableMMS} and \ref{fig:tableMMS_reg2}) but additionally one can choose $u_{\min} = \varepsilon>0$ arbitrarily small. In Figure \ref{fig:controlv2_different_umin} we represent $\|E(t),M(t),F(t)\|$ and $u(t)$ as a function of time for $\varepsilon=10^{-2}$, $1$ and $5$. We see that the curves of $\|E(t),M(t),F(t)\|$ are very similar, the main difference being that $u(t)$ takes more time to converges to 0 as $u_\text{min}$ grows larger. Interestingly, taking $\varepsilon=0$ does not lead to the converges of the equilibrium, suggesting that there is a mathematical bifurcation. Note that being able to take $\varepsilon>0$ arbitrarily small is a much more powerful property than the one given by \eqref{eq:Ustar}.

\begin{table}
    \centering
    \begin{tabular}{|l |c| c| c| c|}
        \hline
         & 200 days & 400 days & 600 days & 800 days \\
        \hline
        average $|E|+|M|+|F|$ & 48,806.91 & 688.68 & 2.47 & 0.002 \\
        \hline
        variance $|E|+|M|+|F|$ & 75,826,146 & 78,547.67 & 1.31 & \(1.61 \times 10^{-6}\) \\
        \hline
        maximum $|E|+|M|+|F|$ & 59,079.04 & 1,130.92 & 4.37 & 0.006 \\
        \hline
        average $|M_{s}|$ & 2,500,000 & 129,308.46 & 2,204.39 & 41.67 \\
        \hline
        variance $|M_{s}|$ & \(3.07 \times 10^{-11}\) & 2,949,520,902 & 5,197,430.54 & \(2.17 \times 10^{-7}\) \\
        \hline
        maximum $|M_{s}|$ & 2,500,000 & 248,387.02 & 10,757.19 & 41.67 \\
        \hline
    \end{tabular}
    \caption{Statistics over 100 simulations with random initial conditions in $[0, 10K]^4$ using control $v_\text{reg}$ (see \eqref{eq:controlv2}) with $u_\text{min}=5$ and $u_\text{max}=300000$ over a duration of 600 days.}
    \label{fig:tableMMS_reg2}
\end{table}

\end{document}